\RequirePackage{fix-cm}
\documentclass{amsart}

\usepackage{graphicx}
\usepackage{mathrsfs}
\usepackage{amsfonts,amsmath}
\usepackage{color}
\usepackage{mathrsfs}

\usepackage{color}
\usepackage{hyperref}

\renewcommand{\P}{{\mathbb P}}

\newcommand{\C}{{\mathbb C}}
\newcommand{\R}{{\mathbb R}}

\newcommand{\I}{{\mathbb I}}

\newcommand \A[1]{{\bf (#1)}}

\def\x{{\mathbf{x}}}

\newtheorem{THM}{Theorem}

\newtheorem{lem}[THM]{Lemma}

\def\0{{\mathbf{0}}}

\newcommand{\E}{{\mathbb E}}

\newcommand{\F}{{\mathcal F}}

\newcommand{\leftB}{{{[\![}}}
\newcommand{\rightB}{{{]\!]}}}
\renewcommand{\S}{{\mathbb S}}

\theoremstyle{plain}
\newtheorem{theorem}{Theorem}[section]
\newtheorem{lemma}[theorem]{Lemma}

\theoremstyle{definition}
\newtheorem{remark}{Remark}[section]

\begin{document}

\title{ $L^p $ Estimates for  Degenerate Non-local Kolmogorov Operators
}


\title{Weak well-posedness of multidimensional stable driven SDEs in the critical case
}

\author{Paul-\'Eric Chaudru de Raynal
}

\address{Univ. Grenoble Alpes, Univ. Savoie Mont Blanc,
 CNRS, LAMA, 73000 Chamb\'ery, France\\
pe.deraynal@univ-smb.fr}

\author{St\'ephane Menozzi}

\address{Laboratoire de Mo\'elisation Math\'ematique d'Evry, UMR CNRS 8071,\\
Universit\'e d'Evry Val d'Essonne, Paris-Saclay,
23 Boulevard de France, 91037 Evry, France,\\
and\\
Laboratory of Stochastic Analysis, Higher School of Economics,\\
Pokrovsky Boulevard, 11, Moscow, Russian Federation\\
stephane.menozzi@univ-evry.fr}

\author{Enrico Priola}
\address{Universit\`a di Pavia, Dipartimento di matematica, Via Adolfo Ferrata 5, 27100 Pavia. enrico.priola@unipv.it}

\maketitle

\begin{abstract} 
  We establish weak well-posedness for critical symmetric stable driven SDEs in $\R^d$ with additive noise $Z$, $d \ge 1$.  Namely, we study the case where the stable index of the driving process $Z$ is  $\alpha=1 $ which exactly corresponds to the order of the drift term having the coefficient $b$ which is continuous and bounded. 
In particular, we cover the cylindrical case when 
 $Z_t= (Z^1_t, \ldots, Z^d_t)$ and $Z^1, \ldots, Z^d$ are independent one dimensional Cauchy processes. 
 Our approach relies on $L^p$-estimates for stable operators and uses perturbative arguments.  
\end{abstract}

\vspace*{1cm}
\\
\hspace*{.3cm}

\section{Statement of the problem and main results}
We are interested in proving well-posedness  for the martingale problem associated with the following SDE:
\begin{equation}
\label{SDE}
X_t=x+\int_0^t b(X_s)ds + Z_t,
\end{equation}
where $(Z_s)_{s\ge 0} $ stands for a symmetric 
$d$-dimensional stable process of order $\alpha=1 $ defined on some filtered probability space $(\Omega,\F,(\F_t)_{t\ge 0},\P) $ (cf. \cite{bass2004} and the references therein) under the sole assumptions of continuity and boundedness on the vector valued coefficient $b$:
\begin{center}
\A{C}\ The drift $b: \R^d \to \R^d $ is \textit{continuous} and bounded.\footnote{The boundedness of $b$ is here assumed for technical simplicity. Our methodology could apply, up to suitable localization arguments, to a drift $b$ having linear growth.}
\end{center} 
Above, the generator $L $ of $Z $ writes: 
   \begin{eqnarray} 
\label{GEN_DRIVING_Z}  
L\varphi(x)&=&{\rm p. v.}\int_{\R^d\backslash\{0\}
} [\varphi(x+z)-\varphi(x)] \nu (dz), \;\;\; x \in \R^d, \;\;\; \varphi \in C_b^2(\R^d), \notag\\ 
\nu(dz)&=&\frac{d\rho}{\rho^{2}}\tilde \mu(d\theta),\ z =\rho \theta, (\rho,\theta)\in \R_+^*\times {\mathbb S}^{d-1}.  
\end{eqnarray} 
 (here $\langle \cdot, \cdot \rangle$ (or $\cdot$)   and $|\cdot|$ denote respectively the inner product and the norm in $\R^d$).
 In the above equation, $\nu$ is the L\'evy  intensity measure of $Z,$ ${\mathbb S}^{d-1}$ is the unit sphere of $\R^d $ and $\tilde \mu $ is a spherical measure on $\mathbb S^{d-1} $. It is well know, see e.g. \cite{sato:99} that the L\'evy exponent $\Phi $ of $Z$ writes as:
\begin{equation}
\label{LEVY_EXP}
\Phi(\lambda)=\E[\exp(i\langle \lambda, Z_1\rangle)]=\exp\Big(-\int_{\mathbb S^{d-1}}|\langle \lambda,\theta\rangle| \mu(d\theta)\Big),\;\; \lambda \in \R^d,
\end{equation}
where $\mu= c_1 \tilde \mu $, for a positive constant $c_1$, is the so-called spectral measure of $Z $.
We will assume some non-degeneracy conditions on $  \mu$. Namely we introduce assumption 
\begin{trivlist}
\item[\A{ND}] There exists $\kappa\ge 1 $ s.t. 
\begin{equation}
\label{ND_MES}
\forall \lambda\in \R^d,\  \kappa^{-1}|\lambda|\le \int_{\mathbb S^{d-1}}|\langle \lambda,\theta\rangle| \mu(d\theta)\le \kappa|\lambda|.
\end{equation}
\end{trivlist}
Notably, no  regularity on the spectral measure itself is assumed.
 In particular, we do not assume that $\mu$ has a density with respect to the Lebesgue measure on  $\mathbb S^{d-1} $. When such density exists and it  is constant, we get, up to multiplicative factor,  the usual fractional Laplacian $L= \Delta^{\frac 12} $. On the other hand, we can also consider  the singular measure $\mu = \sum_{i=1}^d \frac 12 (\delta_{e_i}+\delta_{-e_i}) $, where  $(e_i)_{i\in \leftB 1,d\rightB} $ denotes  the canonical basis of $\R^d $, which corresponds  to the cylindrical fractional Laplacian $\sum_{i=1}^d (\partial_{x_i x_i}^2)^{\frac 12} $
  (cf. \cite{bass:chen:06}).

In dimension one,  \eqref{SDE} is investigated  in the seminal paper \cite{tana:tsuc:wata:74}; the authors  prove that \eqref{SDE} is  well-posed 
if  $b: \R \to \R$ is continuous and bounded and the  exponent $\Phi : \R^d \to \C$ of the L\'evy process $Z= (Z_t)$ verifies $(\text{Re}\Phi (\lambda))^{-1} = O (|\lambda|^{-1})$ as $|\lambda| \to \infty$. Note that, still in dimension one, uniqueness in law implies pathwise uniqueness.
For $d \ge 1$ assuming that  $\mu$ has a density,  well-posedness of \eqref{SDE} follows by the results in \cite{komatsu}.\\

\def\ciao1{  
Even though the approach presented below would also work in dimension one, we assume here that $d> 1 $, i.e. that we are in the strictly multi-dimensional case. This is mainly motivated by the fact that weak-well-posedness in the scalar critical case, as well as some related counter-examples in the super-critical case involving stability indexes lower than one, have been investigated in the famous work of Tanaka \textit{et al.}
\cite{tana:tsuc:wata:74}.
}

   Now, the operator which is at the moment only formally associated with the dynamics in \eqref{SDE} writes for all $\varphi \in C_b^2(\R^d) $, and $x\in \R^d$:
\begin{eqnarray}\label{GEN_SDE}
{\mathcal L} \varphi(x)&=&{\rm p. v.}\int_{\R^d\backslash\{0\}
} [\varphi(x+z)-\varphi(x)] \nu (dz)+\langle b(x),D\varphi(x)\rangle\notag\\
&=& L\varphi(x)+\langle b(x),D\varphi(x)\rangle,
  \end{eqnarray} 
The current setting is said to be critical because, roughly speaking, the two terms in \eqref{GEN_SDE} have the same order one. Therefore it is not clear that the smoothing properties of the semi-group generated by the 
non-local 
operator are sufficient to regularize a transport term. 

If the drift $b$ is itself H\"older continuous and bounded, stronger assumption than  \A{C}, then the well-posedness of the martingale problem for ${\mathcal L} $ can be  established  following \cite{miku:prag:14} (see also Section 3 in \cite{prio:12}).   
  For unbounded H\"older drifts this property follows from the Schauder estimates established in \cite{chau:meno:prio:19}. All these results are based on perturbative techniques  which exploit that the singularities of the corresponding heat-kernel serving as a proxy appearing in the analysis can be absorbed thanks to the H\"older continuity.
 However,  under the sole continuity condition in \A{C} those singularities can only be \textit{averaged}. This is why in the current framework we will exploit $L^p$ estimates for singular kernels. 
 We can for instance mention those of \cite{huan:meno:prio:19} in the more general degenerate Kolmogorov setting.  We also remark  that the  proof of Lemma \ref{KRYLOV_POUR_TOUS} (a Krylov's type estimate on the  resolvent) seems to be of independent interest. 
 
Results on weak well-posedness when $Z$ is  non-degenerate symmetric stable  with  $\alpha >1$ and $b$ is in some $L^p$-spaces are available (see \cite{jin} and the references therein). We also mention \cite{zhao} who investigates  weak well-posedness when $\alpha <1$ for a subclass of  symmetric stable processes $Z$ assuming that the $(1- \alpha)$-H\"older norm of $b$ is   small enough.

\def\ciao2{
Finally, we give a brief overview of existence results for Equation (1) with measurable coefficients a and b  known for some particular cases. 

In dimension one \cite{zanzotto} investigates \eqref{SDE} without drift (i.e., $b=0$) It gives a necessary and sufficient conditions for the existence of the weak solution (see also the references therein).

We  mention  a related existence result in \cite{kurenok} which 
states that the SDE
$$
dX_t = b(t, X_t) dt + dZ_t, \;\;\; X_0 = x
$$
has a weak solution if $b$ is a bounded measurable function and the  exponent $\Phi : \R^d \to \C$ of  $Z= (Z_t)$ verifies  
 $(\text{Re}\Phi (\lambda))^{-1} = o (|\lambda|^{-1})$ as $|\lambda| \to \infty$. Weak uniqueness under such assumptions is still not clear.
}  

Finally we  note that    equivalence between 
weak solutions to a class of SDEs driven by Poisson random measures and solutions to the  martingale
problem for a class of  non-local operators is investigated in \cite{kurtz2011}.

\medskip  
  Let ${\mathscr D}(\R_+,\R^d) $ be the Skorokhod space of all c\`adl\`ag functions from $\R_+$ into $\R^d $ endowed  with the Skorokhod topology and consider the canonical process  $(X_t)$, $X_t(\omega)$ $ = \omega(t)$, for any  $\omega \in 
 {\mathscr D}(\R_+,\R^d)$. 
  Let ${\mathcal L}$ be defined as in \eqref{GEN_SDE}. 
 
 Let us fix a Borel probability measure $\mu$ over $\R^d$.  A solution to the $(  {\mathcal L}, \mu)$-martingale problem  is  a probability measure $\P = \P^{\mu}$ on the Skorokhod space such that $\P (X_0 \in B) = \mu (B)$, $B \in {\mathcal B}(\R^d)$ and, moreover,  
 for any function  $\varphi\in C_b^2(\R^d) $,
$$
M_t=\varphi(X_t)-\varphi(X_0)-\int_0^t {\mathcal L} \varphi(X_s)  ds, \ t\ge 0, 
$$
is a $\P^{\mu} $-martingale (with respect to the canonical filtration). 
 For a comprehensive study of such martingale problems see Chapter 4 in \cite{ethi:kurz:97}. 

  The martingale problem for $\mathcal L$ is well-posed, if  
   for any initial distribution $\mu$ there exists a unique (in the sense of finite-dimensional distributions) solution to the $(  {\mathcal L}, \mu)$-martingale problem.   
Our main result is the following one.
\begin{THM}
\label{THM_WEU}
Under \A{ND} and \A{C}  the martingale problem for ${\mathcal L}$ (cf.  \eqref{GEN_SDE}) is well-posed. 
  \end{THM}
 
\section{Main Steps for the proof of Theorem \ref{THM_WEU}}
Our approach is the following. Let us first assume that $b$ does not vary much, i.e.
\begin{trivlist}
\item[\A{C${}_\varepsilon$}] Additionally to \A{C}, we assume that there exists $\varepsilon \in (0,1) $, $b_0\in \R^d$ s.t.
\begin{equation}\label{VARIE_PEU}
|b(x)-b_0 | \le \varepsilon,\;\;\; x \in \R^d.
\end{equation}
\end{trivlist} 
The previous hypothesis means that we can choose $b_0 \in \R^d$ and $\varepsilon >0$ small enough such that \eqref{VARIE_PEU} holds.
 In this case, let us introduce, for any given starting point $x\in  \R^d $, the following  \textit{frozen proxy} process:
\begin{equation}
\label{PROXY_PROC} 
{{\tilde X}^x}_t=x+ b_0 t + Z_t,\ t\ge 0,
\end{equation}
where we recall that $(Z_s)_{s\ge 0} $ stands for a symmetric non-degenerate  $d$-dimensional stable process of order $\alpha=1 $ defined on some filtered probability space $(\Omega,\F,(\F_t)_{t\ge 0},\P) $. The generator of $(\tilde X_t^x)$ 
   writes for all $\varphi\in C_b^2(\R^d) $ as:
\begin{eqnarray}
\label{GEN_PROXY}  {\mathcal L}^{b_0} \varphi(x) 
= \tilde {\mathcal L}\varphi(x)&=&{\rm p. v.}\int_{\R^d\backslash\{0\}
} [\varphi(x+z)-\varphi(x)] \nu (dz)+\langle b_0 ,D\varphi(x)\rangle\notag\\
&=& L\varphi(x)+\langle b_0 ,D\varphi(x)\rangle, \;\; x\in \R^d,  
\end{eqnarray}
using again the notation \eqref{GEN_DRIVING_Z}.

Using the regularity properties of the density of $Z_t$, 
 under \A{ND}, we find easily that  ${\tilde X}^x_t $ admits for $t>0 $ a $C^{\infty}$-density $\tilde p(t,x,\cdot) $ which can be expressed as:
\begin{equation}
\label{CORRESP_DENS}
\tilde p(t,x,y)=p_{Z_t}\big(y-x- b_0 t\big),\;\;  t>0.
\end{equation}
Following the proof of  Lemma 4.3 in \cite{huan:meno:prio:19} (see also Section 4.2   in \cite{chau:meno:prio:19}) we can as well show the following result (the sketch of the proof is postponed to Appendix).
\begin{lemma}[Controls on the frozen density]\label{LEMME_DENS}
The density $p_{Z_t}$ and its derivatives satisfy the following integrability properties. There exists a constant $C_1:=C_1(\A{ND})\ge 1 $ s.t. for all multi-index  \textbf {$   \beta$}, $|{\bf \beta}|\le 2 $, 
\begin{equation}
\label{DER_DENS_BOUND}
  |\partial_{z}^{\bf \beta} p_{Z_t}(z)|\le \frac{C_1}{t^{| {\bf \beta }|}}  \bar q(t,z), \;\;\;  z \in \R^d,
\end{equation}
where $\bar q(t,\cdot) $ is a probability density s.t. $t^{-d}\,  \bar q(1, t^{-1} x) = \bar q(t,x)$, where $\bar q (1, \cdot) \in L^p(\R^d)$, $p \in [1, \infty)$,  and 
for all $\gamma<1 $, there exists $C_\gamma$ s.t. 
\begin{equation}
\label{INT_BOUND} 
\int_{\R^d}|z|^\gamma \bar q(t,z)\le C_\gamma t^\gamma,\;\; t>0.
\end{equation}
Moreover, the density $\bar q $ enjoys the following property (cf. formula (4.26) in the proof of Lemma 4.3 in \cite{huan:meno:prio:19}):  fix $K\ge 1$; there exists $C:=C(K)$ s.t. for all $t>0,x,y\in \R^d$ with $|y|/t\le K $ then:
\begin{equation}   
\label{PERTURB_DIAG_DENS_BD}
\bar q(t,x+y)\le C\bar q(t,x).
\end{equation}
Also, for all  $t>0,\ z\in \R^d $, $c_0 >0,$ $h\in \R^d$ with $|h|\le c_0 t $, 
\begin{equation}\label{LAP_FRAZ}
|\partial_{z}^\beta p_{Z_t}(z+h)|\le \frac{\tilde C}{t^\beta}\bar q(t,z),\ \tilde C:=\tilde C(C_1,c_0),\ |\Delta^{\frac 12}  p_{Z_t}(z)|\le \frac{C_1}t\bar q\big(t,z\big), 
\end{equation}
and for all $\beta\in (0,1) $, there exists $C_\beta\ge 1 $ s.t. for all $(z,z')\in (\R^d)^2 $:
\begin{equation}
\label{DIFF_LAP_FRAZ}
|\Delta^{\frac 12} p_{Z_t}(z)-\Delta^{\frac 12} p_{Z_t}(z')|\le \frac{C_\beta}t\left( \frac{|z-z'|}{t}\right)^\beta\Big(\bar q\big(t,z\big)+\bar q\big(t,z'\big)\Big),
\end{equation}
as well as 
\begin{equation}
\label{DIFF_LAP_GRAD}
|D p_{Z_t}(z)-D p_{Z_t}(z')|\le \frac{C_\beta}t\left( \frac{|z-z'|}{t}\right)^\beta\Big(\bar q\big(t,z\big)+\bar q\big(t,z'\big)\Big).
\end{equation}
\end{lemma}

For a fixed parameter  $\lambda>0$,
we now introduce the resolvent associated with the frozen proxy process in \eqref{PROXY_PROC}.  Namely, for all $x \in \R^d $ and $f \in C_0^{\infty}(\R^d) $:
\begin{equation}
\label{GREEN}
\tilde R^\lambda \varphi (x):=\int_{0}^{+\infty} \exp(-\lambda t)\E[ ({\tilde X}^x_t)] dt=\int_{0}^{+\infty}\!\!\! dt  \exp(-\lambda t)\int_{\R^d}\tilde p(t,x,y) \varphi(y) dy.
\end{equation}
It is clear that the above function is the unique  classical solution of the  PDE:
\begin{equation}
\label{FROZEN_PDE}
 \tilde {\mathcal L}u(x)-\lambda u=-f(x),\ x\in  \R^d.
\end{equation} 
It also satisfies, for the smooth source considered, some Schauder estimates (see e.g. \cite{prio:12}). Eventually, we also derive from Lemma \ref{LEMME_DENS} the following important pointwise estimate. For all $ p>d$ there exists $C_p$ s.t. for all $x\in \R^d $,
\begin{equation}
\label{CTR_GREEN_LP_RIS}
|\tilde R^\lambda f (x)|\le C_p (1+\lambda^{-1})\|f\|_{L^p}, \;\;  f \in C_0^{\infty}(\R^d).
\end{equation}
Indeed, from \eqref{DER_DENS_BOUND}, one gets (using the H\"older inequality):
\begin{eqnarray*}
|\tilde  R^\lambda f(x)|&\le& \int_{0}^{+\infty} dt \exp(-\lambda t)\int_{\R^d} |f(y)| \frac{C}{t^d} \bar q(1,\frac y t) dy\le C\int_0^{+\infty} \frac{dt}{t^{\frac d p}} \exp(-\lambda t)\|f\|_{L^p}\\
&\le & C_p(1+\lambda^{-1})\|f\|_{L^p}.
\end{eqnarray*}

 We are actually interested in   the corresponding PDE with variable coefficients involving  $\mathcal L$ which writes:
\begin{equation}
\label{VARIABLE_PDE}
\mathcal Lu(x)-\lambda u(x)=- f (x),\ x\in \R^d.
\end{equation}
The idea  is  to express a solution of \eqref{VARIABLE_PDE} in terms of the solutions of \eqref{FROZEN_PDE} which are well understood. 
To this end we first write:
\begin{equation}
\label{VARIABLE_PDE_WITH_GREEN}
 {\mathcal L} \tilde R^\lambda  f(x)-\lambda \tilde R^\lambda  f (x) =-f(x)+{\mathcal R}  f (x)=:-(I-{\mathcal R}) f (x),\ x\in \R^d,
\end{equation}
where in the above equation the remainder operator ${\mathcal R}$ writes for all $ f\in C_0^{\infty}(\R^d)$:  
\begin{equation}
\label{DEF_REMAINDER}
\mathcal Rf(x)=({\mathcal L}-\tilde {\mathcal L})\tilde R^\lambda f(x) =\langle b(x) - b_0  ,D\tilde R^\lambda f(x)\rangle.
\end{equation}
The point is that a \textit{formal} solution of \eqref{VARIABLE_PDE} is provided by the expression $\tilde R^\lambda \circ(I-\mathcal R)^{-1}  f (x) $ if $I-\mathcal R $ can be inverted on a suitable function space (see Section 4.1 for more details). A sufficient condition is that the remainder operator $\mathcal R$ has sufficiently small associated norm. 
Observe that for all $t>0 $,
$$|D\int_{\R^d} \tilde p(t,x,y)f(y)dy|\le Ct^{-1}\| f\|_{\infty}.
$$ 
Under our sole continuity conditions \A{C${}_\varepsilon $}, we cannot 
absorb such a time singularity pointwise (this could be done  in the H\"older framework).

We are thus naturally in the framework of singular integrals, i.e. explosive contribution in time, for which one can expect the time singularity to be absorbed through averaging. The associated natural function spaces to be considered are thus the Lebesgue spaces $L^p( \R^d)$.  

By Lemma \ref{LEMME_DENS} 
we will  actually derive the following theorem 
 which is  the main technical tool to establish the well-posedness for \eqref{SDE}. 
\begin{THM}[$L^p $-estimates for the resolvent] \label{MLP}Under \A{ND} we have, $\lambda >0,$ 
\begin{equation}
\label{BOUND_LP_GRAD}\|D\tilde R^\lambda f\|_{L^p( \R^d)}\le (1 + \lambda^{-1}) C(p, \A{ND},d, |b_0|)\, \|f\|_{L^p( \R^d)}, \; f\in C_0^{\infty}(\R^d).
\end{equation} 
\end{THM}
We deduce  
\begin{THM}[$L^p $-estimates for the remainder in the critical case]\label{THM_LP}
Assume \A{ND} and \A{C${}_\varepsilon $} hold. Let $\lambda \ge 1$. 
For any $ p\in (1,+\infty)$, there exists $C_p := C(p, \A{ND} , d, \| b\|_{\infty})$  s.t. 
\begin{equation} 
\label{LP_EST}
\|\mathcal Rf\|_{L^p( \R^d)}\le C_p\, \varepsilon \,  \|f\|_{L^p( \R^d)},\;\; \; f\in C_0^{\infty}(\R^d).
\end{equation}
\end{THM}  
Importantly,  assuming \A{C${}_\varepsilon $},  we have: $|b_0| \le |b_0 - b(x)| + |b(x)|$ $\le 1 + \| b\|_{\infty}$;  hence 
 the   estimate in \eqref{LP_EST}  is independent of   $b_0  $.  
 
 Theorem \ref{THM_LP}  gives that   $(I - \mathcal R)$ is invertible on $L^p( \R^d) $ if $\varepsilon $ is sufficiently small.
\begin{proof} 
Write for all $x\in  \R^d$:
\begin{eqnarray*}
|{\mathcal R}f(x)|\le |b_0 -b(x)||D\tilde R^\lambda f(x)| \underset{\eqref{VARIE_PEU}}{\le} \varepsilon |D\tilde R^\lambda f(x)|.
\end{eqnarray*}
Hence, 
\begin{equation}\label{PREAL_LP_BOUND}
\|\mathcal Rf\|_{L^p( \R^d)}{\le} \varepsilon \|D\tilde R^\lambda f\|_{L^p( \R^d)}.
\end{equation}
 Plugging \eqref{BOUND_LP_GRAD} in \eqref{PREAL_LP_BOUND} yields the result.
\end{proof}

The previous estimates then allow to establish that the martingale problem associated with $ {\mathcal L}$  is well-posed under the assumption \A{C${}_\varepsilon $}. From this first result, we can get rid of the \textit{almost constant} coefficients through the continuity assumption and a localization argument. These points are detailed in Section \ref{SEC_WP} below.
 
\section{Proof of the main $L^p$ estimates of Theorem \ref{MLP}} 

To get  the $L^p$-estimates in \eqref{BOUND_LP_GRAD}, we will adopt the Coifmann and Weiss approach \cite{coif:weis:71} as in \cite{huan:meno:prio:19}. It therefore suffices to establish the two following lemmas. 
\begin{lemma}[Global $L^2$-estimate]
\label{LEMME_L2}
 There exists a positive constant $C_2:=C_2(\A{ND})$ such that, for all $\lambda>0 $  and for all $f \in L^2( \R^{d})$,
$$
\|D \tilde R^\lambda  f\|_{L^2(\R^d)}
 \le C_2 \|f\|_{L^2(\R^d)}.
$$
\end{lemma}
We mention that this estimate  would hold under weaker assumptions than \A{ND}. In particular, no symmetry would a-priori be needed. 
\begin{lemma}[Deviation Controls]
\label{LEMME_DEV} There exist constants $K$ and $C$ possibly depending on $\A{ND}$ and $|b_0|$,   s.t. for all 
$\xi, x \in \R^d$  the following control hold:
\begin{align*}
\int_0^{+\infty} dt \exp(-\lambda t) \int_{  |x-y|\ge K |x-\xi| }|D\tilde p (t,x,y)-D\tilde p (t,\xi,y) | dy  \le C(1+\lambda^{-1}).
\end{align*}
\end{lemma}
 Indeed, up to a direct symmetrization of the singular kernel involved, since the above estimates still hold for the adjoint operator, we readily derive that \eqref{BOUND_LP_GRAD} follows from the controls of Lemmas \ref{LEMME_L2}, \ref{LEMME_DEV} and Theorem 2.4 in  Chapter III of \cite{coif:weis:71}. To this purpose we also need an additional truncation procedure to separate the singular and non-singular part of the kernel (see e.g. Lemma 3.2 and eq. (3.11) in \cite{huan:meno:prio:19}).  
 
\begin{remark}{
We  point out that, in the critical case { the 
$L^p$-estimates for the fractional Laplacian or the gradient applied to the resolvent are the same.}  Indeed, in the Fourier space their symbol are respectively $|\xi| $ and $-i\xi $. The $L^2$-estimate will not make any difference and as far as the deviations are concerned, we recall from \eqref{DIFF_LAP_GRAD}, \eqref{DIFF_LAP_FRAZ} and \eqref{CORRESP_DENS} that both singular kernels share the same  density estimate. }
\end{remark}  
\def\ciao
{
As a direct consequence of the previous estimate, for a non degenerate and bounded scalar diffusion coefficient $\sigma(w)$, setting $L_wf(x)={\rm p.v.}\int_{\R^d}\big(f(x+\sigma(w)z)-f(x)\big) \frac{dz}{|z|^{d+1}} $, the same type of estimate holds i.e. $\|(L_x-L_{x_0})\tilde R^\lambda f \|_{L^p}\le C |\sigma(x)-\sigma(x_0)| \|f\|_{L^p}$. Indeed, in such a case, it is readily seen that $L_wf(x)=|\sigma(w)| \Delta^{\frac 12} f(x) $. Thus, assumption \A{C${}_\varepsilon $} and the previous $L^p$-estimates give the statement. In the multidimensional setting the proof is more involved we can e.g. refer to theorem 2.2 in \cite{komatsu} in the absolutely continuous case.
Hence, in the following we only prove the results in Lemmas \ref{LEMME_L2}, \ref{LEMME_DEV} for the fractional Laplacian.  
} 
\subsection{Proof of Lemma \ref{LEMME_L2}}
Introduce first for $\eta>0 $, 
\begin{equation}\label{RES_WITH_EPS}
\tilde R_\eta^\lambda f(x)=\int_{\eta}^{+\infty} dt \exp(-\lambda t)\int_{\R^d}\tilde p(t,x,y)f(y) dy.
\end{equation} 
We start from the representation of the density $\tilde p $ obtained under \A{ND} through Fourier inversion from \eqref{CORRESP_DENS}. Namely, for all $t>0, (x,y)\in (\R^d)^2$,
\begin{eqnarray*}
 \tilde p(t,x,y)&=&
 \frac{1}{(2\pi)^d}\int_{\R^{d}} \exp\Big(-i\langle (y-(x+b_0 t)),p\rangle \Big)
  \cdot  \exp\Big(-t\int_{\S^{d-1}}|\langle  p, \xi\rangle| \mu(d\xi)\Big)dp.
 \end{eqnarray*}
Let $k=1, \ldots, d$ and $D_k = \partial_{x_k}$. We can compute, for $f\in{\mathscr S}(\R^{d})$ (Schwartz class of $\R^{d}$), the Fourier transform:
 \begin{eqnarray*} 
 \R^d \ni \zeta \mapsto  \mathcal{F}( D_k \tilde R_\eta^\lambda f) (\zeta) = \int_{\R^{d}}  e^{i \langle \zeta, x \rangle} D_k \tilde R_\eta^\lambda f(x) d x = (-i\zeta_k) \mathcal{F}(\tilde R_\eta^\lambda f)(\zeta)
  \\ 
=(-i\zeta_k) \int_{\R^n}e^{-i\langle \zeta, x\rangle}\Big(\int_{\eta}^{+\infty}\exp(-\lambda t)\int_{\R^d}\tilde p(t,x,y)f(y)dy dt \Big)dx.
\end{eqnarray*}
Using the Fubini theorem and \eqref{CORRESP_DENS}, we derive: 
\begin{eqnarray*}
\mathcal{F}(D_k  \tilde R_\eta^\lambda f)  (\zeta)
&=&(-i\zeta_k)  \int_{\eta}^{+\infty}\exp(-\lambda t ) \int_{\R^d}\exp(-i\langle \zeta, y) f(y) \\
&&\times\Big( \int_{\R^d}\exp(-i\langle \zeta, x-y\rangle)
 p_{Z_t} \big(y-(x+b_0 t) \big)dx \Big) dy dt\\
&=&(-i\zeta_k) \Big(\int_{\eta}^{+\infty}\exp(-\lambda t ) \int_{\R^d}\exp(-i\langle \zeta, y) f(y) \\
&&\exp(i t\langle \zeta,b_0  \rangle)\times\Big( \int_{\R^d}\exp(-i\langle \zeta, \tilde x\rangle)
p_{Z_t} (-\tilde x)d\tilde x \Big) dy dt,
\end{eqnarray*}
setting $\tilde x=x+b_0 t-y $ for the last identity. Recalling that by symmetry $p_{ Z_t }(-\tilde x)=p_{ Z_t}(\tilde x) $, we finally get
\begin{eqnarray*}
&&\mathcal{F}(D_k  \tilde R_\eta^\lambda f)  (\zeta)\\
&=& (-i\zeta_k)  \int_{\eta}^{+\infty} \exp\big(t(-\lambda+i\langle \zeta,b_0 \rangle)\big) \mathcal{F}(f)( \zeta)\mathcal{\F}(p_{Z_t})(t,\zeta)dt\\
&=& (-i\zeta_k)  \int_{\eta}^{+\infty} \exp\big(t(-\lambda+i\langle \zeta,b_0 \rangle)\big)  \mathcal{F}(f)(s, \zeta)\exp\left(- t\int_{\S^{d-1}} |\langle \zeta,\xi \rangle|\mu(d\xi) \right) dt,
\end{eqnarray*}
where ($\mathcal{F}(f)(\cdot), \mathcal{F}(p_{ Z_t  })(\cdot)$ denote the Fourier transforms of $f(\cdot),p_{Z_t }(\cdot)$.

Let us now com\-pu\-te $\|\mathcal{F}(D_k \tilde R_\eta^\lambda f)\|_{L^2(\R^d)}$. From the non degeneracy of $\mu$ we have:
\begin{equation}
\label{THE_REF_FOR_LATER}
|\mathcal{F}(D_k \tilde R_\eta^\lambda f)  (\zeta) | 
\le C |\zeta| \int_{0}^{+\infty} \exp(-\lambda t)  |\mathcal{F}(f)( \zeta)|\exp\left(-C^{-1}t| \zeta| \right)dt.
\end{equation}
 For the $L^2$-norm of $\mathcal{F} (D \tilde R_\eta^\lambda f)$, using the Cauchy-Schwarz inequality, we obtain:
\begin{eqnarray*}
\|\mathcal{F}(D_k \tilde R_\eta^\lambda f)\|_{L^2(\R^d)}^2
&\le& C   \int_{\R^{d}} \Bigg(   |\zeta|
 \,\int_{0}^{+\infty} \exp(-\lambda t)  |\mathcal{F}(f)( \zeta) |^2 \exp\left(- C^{-1}t| \zeta| \right)dt\Bigg)\\
 &&\times \Bigg(|\zeta|\int_0^{+\infty}\exp(-\lambda  t) \exp\left(- C^{-1}t| \zeta| \right)dt\Bigg)d\zeta\\
 &\le &C   \int_{\R^{d}}    |\zeta|
 \,\int_{0}^{+\infty} \exp(-\lambda t)  |\mathcal{F}(f)( \zeta) |^2 \exp\left(- C^{-1}t| \zeta| \right)dt d\zeta\\
 &\le & C   \int_{\R^{d}} |\mathcal{F}(f)( \zeta) |^2 \Big(\int_{0}^{+\infty} |\zeta|\exp\left(- C^{-1}t| \zeta| \right)dt\Big) d\zeta\\
 &\le& C   \int_{\R^{d}} |\mathcal{F}(f)( \zeta) |^2d\zeta=C\|\mathcal{F}(f)\|_{L^2(\R^d)}.
 \end{eqnarray*}
The assertion now follows for $f\in {\mathscr S}(\R^{d}) $ from the Plancherel isometry, and for $f\in L^2(\R^{d}) $ by a density argument. From the uniformity of the previous controls, the final statement can eventually be derived letting $\eta $ go to zero through weak convergence arguments.  \hfill $\qed $

\subsection{Proof of Lemma \ref{LEMME_DEV}}
 Let us denote $\rho:=|x-y|, \gamma:=|x-\xi| $. 
From \eqref{LAP_FRAZ}, \eqref{DIFF_LAP_FRAZ}  and  the correspondence \eqref{CORRESP_DENS}, we 
%
%
 easily get that, for a fixed $\beta\in (0,1) $:
\begin{eqnarray}
\label{CTR_DEV}
&&|D\tilde p(t,x,y)-D\tilde p(t,\xi,y)|
= |Dp_{Z_t}\big(y-x-b_0 t \big)-Dp_{Z_t}\big(y-\xi-b_0 t \big)|\notag\\
&\le& \frac{C_\beta}{t}\left(\frac{|x-\xi|}{t}\right)^\beta\big(\bar q(t, y-x-b_0 t)+\bar q(t,y-\xi-b_0 t) \big)\notag\\
&\le & \frac{C_\beta}{t}\left(\frac{|x-\xi|}{t}\right)^\beta\big(\bar q(t, y-x)+\bar q(t,y-\xi) \big)
\end{eqnarray}
using as well the first bound in \eqref{LAP_FRAZ} for the last identity. 
Namely, a  {diagonal} perturbation does not affect the density estimate.
  We therefore derive: 
\begin{eqnarray} \label{sii} \nonumber
I&:=& \int_0^{+\infty} dt \exp(-\lambda t) \int_{  |x-y|\ge K |x-\xi| }| D\tilde p (t,x,y)- D \tilde p (t,\xi,y) | dy  
\\
&\le &C \int_0^{+\infty} dt \exp(-\lambda t)\int_{  \rho >K\gamma} \frac{1}{t}\left(\frac{|x-\xi|}{t} \right)^\beta
\Big(\bar q\big(t, x-y \big) +\bar q\big(t, \xi-y \big) \Big)  dy .
\nonumber\\
&=:&I_1+I_2.
\end{eqnarray}
We can rewrite, 
\begin{eqnarray}
I_1&=&C\gamma^\beta \int_{\rho>K\gamma} \!\!\! dy \int_0^{+\infty} dt (\I_{t\ge \gamma}+\I_{t< \gamma})   \frac{dt}{t^{1+\beta}}  \exp(-\lambda t) \bar q(t,x-y)=:I_{11}+I_{12}.\notag\\ \label{SPLIT_I1}
\end{eqnarray}
On the one hand, assuming w.l.o.g. that $\gamma\le 1 $, we get
\begin{eqnarray}
I_{11}&\le &C\gamma^\beta\int_{\gamma}^{+\infty}dt  \frac{\exp(-\lambda t)}{t^{1+\beta}} \int_{\R^d} \bar q(t, y-x)d y\le C \gamma^\beta \Big(\int_{\gamma}^1 \frac{dt}{t^{1+\beta}}+\int_{1}^{\infty}\frac{dt}{\gamma^\beta}\exp(-\lambda t) \Big)\notag\\
&\le & C (1+\lambda^{-1}).\label{BD_I11}
\end{eqnarray}
Write now, from the self-similarity properties of $\bar q $ and taking $0<\beta< \eta <1 $:
\begin{eqnarray}
I_{12}&\le &C\gamma^\beta\int_0^{\gamma}dt  \frac{\exp(-\lambda t)}{t^{1+\beta}} \int_{t|\tilde y|\ge K \gamma} \bar q(1, \tilde y)d \tilde y\notag\\
&\le& C \gamma^\beta \Big(\int_0^{\gamma} \frac{dt}{t^{1+\beta}}\int_{\frac{t |\tilde y|}{\gamma}\ge K} \Big(K^{-1}\frac{t |\tilde y|}{\gamma}\Big)^{\eta} q(1,\tilde y) d\tilde y\Big)\notag\\
&\le& CK^{-\eta}\gamma^{\beta-\eta}\int_0^{\gamma}\frac {dt}{t^{1+\beta-\eta}}\int_{\R^d}|\tilde y|^\eta \bar q(1,\tilde y) d\tilde y\le \tilde C.\label{BD_I12}
\end{eqnarray}
Let us now turn to $I_2$ in \eqref{sii}.
\begin{eqnarray*}
I_2&=& C \gamma^{ \beta} \int_{0}^{+\infty} \frac{dt}{t^{1+\beta}}\exp(-\lambda t)    (\I_{t>  \gamma}+\I_{t\le \gamma})
 \Big(\int_{\R^d} \I_{   
\rho >K\gamma} \, \bar q\big(t, \xi-y \big)  dy \Big)=:I_{21}+I_{22}. \notag\\
\end{eqnarray*}
The contribution $I_{21} $, can be handled as $ I_{11}$ introduced in \eqref{SPLIT_I1} and therefore also satisfies \eqref{BD_I11}.
To analyze $I_{22} $, we first set $z=t^{- 1 }(\xi-y)$ and recall that:
\begin{equation}
K\gamma < \rho=|x-y|\le |\xi-y|+|x-\xi|\le |\xi-y|+\gamma \Rightarrow |\xi-y|\ge (K-1)\gamma,
\end{equation}
from which we deduce that $I_{22}$ can be handled as $I_{12} $ above up to a modifications of the considered constants. Hence, it also satisfies \eqref{BD_I12}. Namely, we eventually have from \eqref{BD_I11} and \eqref{BD_I12} that for all $i\in \{1,2\} $,
$I_i\le C(1+\lambda^{-1})$, which plugged into \eqref{sii} gives the statement.

\section{Well-posedness of the martingale problem}\label{SEC_WP}
We prove in this section our Theorem \ref{THM_WEU}. First under the local condition \A{C${}_\varepsilon $} and then extend it through the continuity assumption through a localization argument.

\smallskip 
{\it Existence} of a solution to the martingale problem for $({\mathcal L} , \mu)$, for any initial distribution $\mu$, can be proved, under  \A{ND} and \A{C${}$},  applying, for instance, Theorem 4.1 in \cite{kuhn}. Such theorem 
 is based on 
 Theorem 4.5.4 in \cite{ethi:kurz:97} about  relations between  the positive maximum principle and the existence of solutions to the martingale problem. 
 Theorem 4.5.4  has been used to prove  existence results for martingale problems associated to pseudodifferential operators with symbols $\Phi(x, \lambda )$ which are continuous in $x$. We only mention here   \cite{hoh} and \cite{kuhn} (see also the references therein).  

{\sl Thus in the sequel we concentrate on the problem of uniqueness. }

\subsection{Uniqueness  of the martingale problem under \A{C${}_\varepsilon $}}


Fix a Borel probability distribution $\mu$ on $\R^d$. Let $\P^\mu$ be any solution to the $(  {\mathcal L}, \mu)$-martingale problem  and denote by $(X_t)_{t\ge 0} $ the associated canonical process on the Skorokhod space. Moreover define $\E^{\P^\mu}= \E^\mu$. 

To derive uniqueness we will crucially rely on the following Krylov type estimate whose proof is postponed to the end of the section.
\begin{lem}[Krylov type bound]\label{LEMMA_KRYL_PER_TUTTI} Let $\P^\mu$ be any solution to the $(  {\mathcal L}, \mu)$-martingale problem. 
Define the corresponding resolvent
\begin{gather} \label{re2}
G(\lambda) f = {\E^\mu}\Big[\int_{0}^{+\infty} \exp(-\lambda s) f(X_s) ds
\Big].
\end{gather}
Then  for all $f\in C_0^{\infty}(\R^d) $, $\lambda \ge 1,$ and for all $p>d $ there exists $C:=C(p, d, \varepsilon, \| b\|_{\infty})  >0 $ s.t. 
\begin{equation}\label{KRYLOV_POUR_TOUS}
| G(\lambda) f | \le C \|f\|_{L^{p}}.
\end{equation}
\end{lem}
Observe that the above bound actually provides by duality that 
 $G(\lambda)$ possesses a density $p^{\lambda} = p^{\lambda, \mu}$ which belongs to $L^q(\R^d)$,
$p^{-1}+q^{-1} =1$.
 Lemma \ref{LEMMA_KRYL_PER_TUTTI} precisely provides a good tool to derive uniqueness. 
 \subsubsection{Derivation of uniqueness} 
We will apply Theorem 4.4.2 in \cite{ethi:kurz:97} (see also Theorem 6.2.3 in \cite{SV}).  This says that if two solutions $\P_1$ and $\P_2$ of the $(  {\mathcal L}, \mu)$-martingale problem  have the same one-dimensional marginal distributions then  they coincide (i.e., they have the same finite-dimensional distributions). 
 
\smallskip   Let $\P^\mu$ be any solution to the $(  {\mathcal L}, \mu)$-martingale problem. Let $\lambda \ge 1$. 
  For $f\in C_0^{\infty}(\R^d) $ consider the resolvent  $\tilde R^\lambda f(x)$ introduced in \eqref{GREEN}. This is a smooth function (see e.g. \cite{prio:12}). Set  $\varphi(x)= \tilde R^\lambda f(x)$.  The idea is now to expand,
$\varphi(X_t):= \tilde R^\lambda f(X_t) $. We know that  
\begin{eqnarray*} 
&&\varphi(X_t)-\varphi(X_0)-\int_0^t {\mathcal L} \varphi(X_s) ds = \tilde R^\lambda f(X_t)-\tilde R^\lambda f(X_0)-\int_0^t    {\mathcal L} \tilde R^\lambda f(X_s) ds
 \end{eqnarray*}
is a $\P^\mu $-martingale. We  write:
\begin{eqnarray*}
&&{{\E^{\mu}}}[ \tilde R^\lambda f(X_t)]-{\E^{\mu}}\tilde R^\lambda f(X_0)={{\E^{\mu}}}\left[\int_0^t    ({\mathcal L} - \lambda ) \,  \tilde R^\lambda f(X_s) ds\right]
 + \lambda {{\E^{\mu}}}\left[\int_0^t      \tilde R^\lambda f(X_s) ds\right]
\\
&=&-{{\E^{\mu}}}\left[\int_0^t  f(X_s) ds\right]
+{{\E^{\mu}}}\left[\int_0^t \Big( {\mathcal L}-\tilde {\mathcal L} \Big)\tilde R^\lambda f(X_s) ds\right] + \lambda {{\E^{\mu}}}\left[\int_0^t      \tilde R^\lambda f(X_s) ds\right].
\end{eqnarray*}
Integrating over $[0, \infty)$ with respect to $e^{- \lambda t} dt$ and using the Fubini theorem we find 
\begin{eqnarray}
\E^{\mu}    \tilde R^\lambda f(X_0)&=& {{\E^{\mu}}}\left[\int_0^{+\infty}  \exp(-\lambda s )f(X_s) ds\right]
-{{\E^{\mu}}}\left[\int_0^{+\infty}  \exp(-\lambda s )\Big( {\mathcal  L}-\tilde {\mathcal L} \Big)\tilde R^\lambda f(X_s) ds\right]\notag\\
&=&{{\E^{\mu}}}\left[\int_0^{+\infty}  \exp(-\lambda s )(I-{\mathcal R})f(X_s) ds\right],\label{ESPRESSIONE_RISOLVENTE}
\end{eqnarray}
where  we have used the remainder operator ${\mathcal R}$ which writes for all $ f\in C_0^{\infty}(\R^d)$:  
\begin{equation}
\label{rem}
\mathcal Rf(x)=({\mathcal L}-\tilde {\mathcal L})\tilde R^\lambda f(x) =\langle b(x) - b_0, D\tilde R^\lambda f(x)\rangle.
\end{equation}
 Let now $\P_1$ and $\P_2$ be two solutions for the $(  {\mathcal L}, \mu)$-martingale problem. Let  $G_i(\lambda) f = {{\E^{\P_i}}}\left[\int_0^{+\infty}  \exp(-\lambda s )f(X_s) ds\right]$, $i=1,2$.
  We find by \eqref{ESPRESSIONE_RISOLVENTE}
 \begin{gather*}
 G_i(\lambda) f   = \E^{\mu}    \tilde R^\lambda f(X_0)
 + 
  G_i(\lambda) {\mathcal R} f.
 \end{gather*}
  Define  $T(\lambda) : C^{\infty}_0 (\R^d)\to \R$,
$
T(\lambda) f = G_1(\lambda) f  - G_2(\lambda) f. 
$
We have  
\begin{align}\label{tl}
    T(\lambda) \,  (I -  {\mathcal R}) f =0,\;\;\; f\in C_0^{\infty}(\R^d).  
\end{align}
By using Lemma \ref{LEMMA_KRYL_PER_TUTTI} 
 we know that $T(\lambda)$  
can be extended to a bounded linear operator from $L^p(\R^d)$ into $\R$ (we still  denote by $T(\lambda)$ such extension).
 By Theorem \ref{THM_LP} we know  that $I-{\mathcal R} $ is invertible on $L^p(\R^d)$ for $p>d$ (under assumption \A{C${}_\varepsilon $}). Hence  choosing $\varepsilon$ small enough (independently of $\lambda \ge 1$) we find $
\| T(\lambda)    \|_{L (L^p(\R^d) ; \R)} =0,$ $  \lambda \ge 1.  
$
 
  We obtain that $  {\E^{\P_1}}f(X_t) =  {\E^{\P_2}} f(X_t)$, $t \ge 0,$ for any $f \in C_0^{\infty}(\R^d)$, by using the uniqueness of the Laplace transform.
  By a standard approximation procedure we also get,
    for any Borel set $B \subset \R^d$, 
\begin{gather*}
\P_1 (X_t \in B) = \P_2 (X_t \in B), \;\;\; t\ge 0. 
\end{gather*}

\subsubsection{Proof of Lemma \ref{LEMMA_KRYL_PER_TUTTI}}

We will adapt  an   approximation technique of resolvents introduced  by N.V. Krylov in the gaussian setting (cf. \cite{krylov85} and Chapter VI in \cite{bass98}).  

Arguing as in Theorem 2.1 of \cite{komatsu73} (see also  Theorem 3.1 in \cite{chenwang} and Section 3 in  \cite{kurtz2011})
 one proves   that on some stochastic
  basis $(\Omega, {\mathcal  F}, ( {\mathcal  F}_t), \P)$ there exists a $d$-dimensional 
 1-stable process $Z$ as in \eqref{SDE}, an ${\mathcal F}_0$-measurable r.v. $X_0$ with law $\mu$ and a    solution $Y = (Y_t)$ 
  to 
   \begin{gather*}
Y_t = X_0 + \int_0^t b(Y_s) ds + Z_t, \;\; t \ge 0,
\end{gather*}
such that the law of $Y$ coincides with  $\P^{\mu}$ on the Skorokhod space. 
 
 Fix $\lambda \ge 1$.
 For any Borel set $C \subset \R^d$ the measure 
\begin{gather*}
\gamma (C) = \E^{}  \int_0^{\infty} e^{- \lambda t } 1_C (Y_t) dt
\end{gather*}
is well-defined. Moreover for any ball $ B(z, r)$ we have
\begin{equation}\label{ee}
\gamma (B(z, r) ) >0,\; z \in \R^d, \, r>0. 
\end{equation}
We argue by contradiction. If for some ball $B= B(z, r)$ we have $\gamma (B)=0$, then for any $T>0$ there exists    $0<t <T$ such that
\begin{gather*}
\E^{}  [1_B (Y_t)] = \P^{}  (Y_t \in B) =0.
\end{gather*}
 We choose $T>0$ such that $T \| b\|_{\infty} < r/2$. There exists  $0<t< T$
 such that   $\P^{}  ( |Y_t - z| < r ) =0$.
We find  
\begin{gather*}
\P^{}  ( |Y_t - z| < r ) \ge  \P^{}  ( |\int_0^t b(Y_s) ds | < r/2 , \,  |Z_t +X_0 - z |<r/2)
\\
= \P^{}  (   |Z_t +X_0 - z | < r/2) >0, 
\end{gather*}
because $X_0$ is independent of $Z_t$ and the support of the distribution of $Z_t$ is $\R^d$ (see, for instance, Theorem 3.4 with $A=0$ in \cite{PSXZ}). We have found a contradiction and so \eqref{ee} holds. 

Recall that  $\E^{\mu}$ denotes expectation with respect to $\P^{\mu}$ and $(X_t)$ is the canonical process. As before 
we consider  the measure $\gamma = \gamma^{\mu,\, \lambda}$: 
\begin{gather*}
\gamma (C) = \E^{\mu}  \int_0^{\infty} e^{- \lambda t } 1_C (X_t) dt,\;\; C \in {\mathcal B}(\R^d). 
\end{gather*}
Now we introduce $\phi (x) = c_d (1  + |x|^{d+1})^{-1}$,  $x \in \R^d$, with $c_d = (\int_{\R^d} \phi(x) dx)^{-1}.$ Note the following bound on the first and second derivatives of $\phi$:
\begin{equation}\label{d35}
 |D \phi(x)| + \|D^2 \phi(x) \|_{\R^d \otimes \R^d} \, \le C_d \, \phi(x),\;\; x \in \R^d.
\end{equation}
For $\delta >0$, we consider $\phi_{\delta}(x) = \frac{1}{\delta^{d}}  \phi (x/ \delta)$, $x \in \R^d$. Using  $\phi_{\delta}(x-y)$ $\ge \frac{1}{2} \, 1_{B(x, \delta)}(y)$ and \eqref{ee}  we can define 
 \begin{equation}\label{w22}
b_{\delta} (x) = \frac{\int_{\R^d} \phi_{\delta} (x-y) b(y) \gamma (dy)}
{ \int_{\R^d} \phi_{\delta} (x-y) \gamma (dy)},\;\; x \in \R^d, \;\; \delta>0. 
\end{equation} 
Using \eqref{d35}   it is straightforward to check that 
$b_{\delta} \in C^{2}_b(\R^d, \R^d),$ i.e. $b_{\delta}$ has first and second bounded and continuous derivatives, $\delta >0.$ Moreover $\| Db_{\delta}\|_{\infty} \le \frac{c}{\delta}$,
 $\| D^2 b_{\delta}\|_{\infty} \le \frac{c}{\delta^2}$.
 
Let now $\varphi \in C^2_b(\R^d)$; by the martingale property:
$$
\E^{\mu} \varphi(X_t)= \E^{\mu} \varphi(X_0)+ \E^{\mu}\int_0^t [{\mathcal L} \varphi(X_s) -\lambda \varphi(X_s)] ds +    \E^{\mu} \lambda \int_0^t\varphi(X_s)] ds.       
$$
Integrating over $[0, \infty)$ with respect to $e^{- \lambda t} dt$ and using the Fubini theorem we find 
$$
  \E^{\mu}\varphi(X_0)= \int_{\R^d} [\lambda \varphi(y)  - {\mathcal L} \varphi(y)]   \gamma(dy).   
$$
Now we replace  $\varphi$ by  the convolution $\varphi * \phi_{\delta}$    so that
(using also \eqref{d3})
\begin{gather*}
\E^{\mu} [\varphi * \phi_{\delta} (X_0)] = \int_{\R^d} [\lambda \varphi * \phi_{\delta}(y)  - {\mathcal L} [\varphi * \phi_{\delta}] (y)]   \gamma(dy)   
\\
= \int_{\R^d} [ \lambda \varphi * \phi_{\delta}(y)  - {L} \varphi * \phi_{\delta} (y) ]  \gamma(dy) - \int_{\R^d} D \varphi(z)  \cdot \int_{\R^d } \phi_{\delta}(z-y) b(y)  \gamma(dy)\,   dz  
\\
= \int_{\R^d} [ \lambda \varphi * \phi_{\delta}(y)  - {L} \varphi * \phi_{\delta} (y) ]  \gamma(dy) - \int_{\R^d} D \varphi(z)  \cdot b_{\delta}(z)  \int_{\R^d}
\phi_{\delta}(z-y)   \gamma(dy)\,   dz  
\\
=  \int_{\R^d} \int_{\R^d} \big( \lambda \varphi(p)  - [ b_{\delta} \cdot D + L]
\varphi (p) \big) \,   \phi_{\delta}(p-y) \gamma(dy) dp.   
\end{gather*} 
Now we consider the operator 
$$
{\mathcal L}_{\delta} =  b_{\delta} \cdot D + L,
$$
and 
 take $\varphi = G_{\delta}(\lambda)f $, $f \in C_0^{\infty}(\R^d)$, where $G_{\delta}(\lambda)$ is  the resolvent of the  martingale solution $\P^{\delta,x}$ starting at $\delta_x$ associated to the operator ${\mathcal L}_{\delta}$:
$$
G_{\delta}(\lambda)f(x) =  \E^{\delta,x}\big[\int_0^{+\infty} 
\exp(-\lambda s)f(X_s) ds \big],\;\; x \in \R^d. 
$$
 Since $b_{\delta} \in C_b^2(\R^d, \R^d)$ it is not difficult to check that 
 that $\lambda G_{\delta}(\lambda)f(p)  - [ b_{\delta} \cdot D + L]
 \, G_{\delta}(\lambda)f (p) =f(p) $, $p \in \R^d $. It follows that 
 \begin{gather*}
\E^{\mu} [ G_{\delta}(\lambda)f  * \phi_{\delta} (X_0)] =  
  \int_{\R^d} f * \phi_{\delta}(y) \gamma(dy).    
\end{gather*}
Hence we get  the crucial  approximation result:
 \begin{gather} \label{dee}
\E^{\mu} [ G_{\delta}(\lambda)f  * \phi_{\delta} (X_0)]   \to \E^{\mu}  \int_0^{\infty} e^{- \lambda t } f (X_t) dt,\;\; \delta \to 0^+,   
\end{gather}
since $(f * \phi_{\delta})$ is uniformly bounded  and converges pointwise to $f$ on $\R^d$.

Now we consider 
 \eqref{ESPRESSIONE_RISOLVENTE}  
 which we write for the solution $\P^{\delta,x}$ to the martingale problem for ${\mathcal L}_{\delta}$ (starting at the delta Dirac in $x$). We write   
  ${\mathcal R}^{\delta}f(x) = \langle b_{\delta}(x) - b_0, D\tilde R^\lambda f(x)\rangle$ and note that 
\begin{equation}\label{w66}
 |b_{\delta}(x) - b_0| < \varepsilon, \;\; x \in \R^d, \; \delta>0.
\end{equation}
Similarly  to \eqref{ESPRESSIONE_RISOLVENTE} we find   
$$
\tilde R^\lambda f(x)
=
\E^{\delta,x} \int_0^{+\infty } \exp(-\lambda s)\int_{\R^d} (I-{\mathcal R}^{\delta}) f(X_s) ds.
$$
 Since $b_{\delta} $ $\in C^2_b(\R^d, \R^d)$ it is known
that the resolvent $G_{\delta}(\lambda)f(x)= \E^{\delta,x}[\int_0^{+\infty} ds \exp(-\lambda s)f(X_s)]$ has a density $p^{\lambda, \delta}(x,y)$, i.e.,
\begin{equation}\label{d22}
 G_{\delta}(\lambda)f(x) = \int_{\R^d} f(y) p^{\lambda, \delta}(x,y) dy,\;\; x \in \R^d.
\end{equation}
Moreover,  $p^{\lambda, \delta }(x,\cdot) \in L^q(\R^d)$, $q > 1$.
  This follows from the estimate
\begin{equation}\label{d3}
 |G_{\delta}(\lambda)f(x)| \le C_{\delta} \| f\|_{L^p(\R^d)}, \;\; f \in C_0^{\infty}(\R^d),\; \delta>0;
\end{equation}   
see Lemma \ref{d66}.   It follows that 
$$
\tilde R^\lambda f(x) = 
 \int_{\R^d} (I-{\mathcal R}^{\delta}) f (y) p^{\lambda, \delta}(x,y) dy. 
$$
 Recall from \eqref{CTR_GREEN_LP_RIS} that for all $p>d$ there exists $C_p$ s.t. for all $\x\in   \R^{d} $ 
\begin{equation*}
\label{CTR_GREEN_PONCT}
|\tilde R^\lambda f(x)|\le C_p(1+\lambda^{-1}) \|f\|_{L^{p}(\R^{d})}, \;\; f\in C_0^{\infty} (\R^{d}).
\end{equation*}
 Moreover, 
 by Theorem \ref{THM_LP} we know  that $I-{{\mathcal R}^{\delta}} $ is invertible on $L^p(\R^d)$ for $p>d$. For any  $g \in L^p(\R^d)$ we find
 \begin{gather*}
 G_{\delta}(\lambda)g(x)= \int_{\R^d} g(y) p^{\lambda, \delta }(x,y) dy
 \\ 
\le  C_p(1+\lambda^{-1})\, \|(I-{{\mathcal R}^{\delta}})^{-1}\, g\|_{L^p}\le \tilde C_p(1+\lambda^{-1})\| g\|_{L^p}.
\end{gather*}  
It follows that, for $\delta >0,$ 
\begin{gather*}
|\E^{\mu} [ G_{\delta}(\lambda)f  * \phi_{\delta} (X_0)] | \le \tilde C_p(1+\lambda^{-1})\| f\|_{L^p},
\end{gather*}
 $f \in C_0^{\infty}(\R^d)$.  Passing to the limit as $\delta \to 0^+$ we get the assertion by \eqref{dee}.
\qed

\subsection{Uniqueness of the martingale problem under \A{C${}$} by a  localization argument} 

We will use Theorem 4.6.2 in \cite{ethi:kurz:97}. 

First recall  that an operator like $\mathcal L$ under assumption \A{C${}_\varepsilon $} has the property that the associated   martingale problem is well-posed for any initial distribution $\mu$. 

By Theorem 4.6.1 in \cite{ethi:kurz:97},  such  operator  under  \A{C${}_\varepsilon $} has the following  additional property:   for any initial distribution $\mu$ and for any open set $U \subset \R^d$ there exists a unique in law probability measure $\P$ on  ${\mathscr D}(\R_+,\R^d) $  
 such that  
   $\P (X_0 \in B) = \mu (B)$, $B \in {\mathcal B}(\R^d)$, 
$$
X(\cdot) = X(\cdot \wedge \tau), \;\;\; \text{$\P$-a.s., where}  
\;\;    
\tau = \tau^U= \inf \{ t \ge 0 \; : \; X_t \not \in U  \}
$$
($\tau = + \infty$ if the set is empty) is a stopping time with respect to the canonical  filtration, and finally   $ \varphi (X_{t \wedge \tau}) - \int_0^{t \wedge \tau}
    {\mathcal L} \varphi (X_s) ds,$ $ t \ge 0,$ $\varphi \in C_b^2(\R^d)$, 
is a martingale with respect to the canonical filtration (recall that $(X_t) = (X(t))$ indicates the canonical process).
 
One says that  under assumption \A{C${}_\varepsilon $} for any open set $U \subset \R^d$, for any initial distribution $\mu$
the  stopped martingale problem for $({\mathcal L}, U, \mu)$ is well-posed.

 {\sl Now only assuming \A{C${}$},   
 we construct a suitable 
 covering $(U_{j})_{j \ge 1}$ of open sets in $\R^d$ such that for any initial distribution $\mu$ the  stopped martingale problem for $({\mathcal L}, U_j, \mu)$ is well-posed, for any $j \ge 1$. According to Theorem 4.6.2 in \cite{ethi:kurz:97} we conclude that the (global) martingale problem for $\mathcal L$ is well-posed.   
}

To construct such covering of  $\R^d$ we note that 
  by the continuity of  $b$
 we can find  a sequence  $(x_j) \subset \R^d$,
$j \ge 1$, and numbers $\delta_j>0$
 such that the open balls $U_j = B(x_j,  \delta_j)$ of center $x_j$ and radius $\delta_j$
  form a covering for $\R^d$ and moreover
  we have  $
  |b(x) - b(x_j)| < \varepsilon$ (cf. \eqref{VARIE_PEU})
 for any $x \in B(x_j, 2 \delta_j)$, $j \ge 1$.     
    
 The balls
 $\{ B(x_j,  \delta_j) \}_{j \ge 1}$ give the required covering $\{ U_j\}_{j \ge 1}$.  
 Let us define   operators ${{\mathcal L}}_j$ such that
 \begin{align} \label{bo1}
 {\mathcal L}_j \varphi (x) = {\mathcal L} \varphi (x),\;\; x \in U_j,\;\; \varphi \in C^{2}_b(\R^d),
 \end{align}
 and such that  each ${\mathcal L}_j$ verifies \A{C${}_\varepsilon $}.  
  We fix $j \ge 1$ and consider
  $\eta_j \in C^{\infty}_0(\R^d)$ with $0 \le \eta_{j} \le 1$,
 $\eta_{j} =1$ in $B(x_j, \delta_j)$ and  $\eta_j =0$ outside $B(x_j, 2
 \delta_j )$. Now define
 $$
b^j(x) :=  \eta_{j}(x) b(x) +  (1 - \eta_{j}(x)) b (x_j).
 $$
 
 It is easy to see that  $b^j(x) = b(x)$, $x \in U_j$ and 
   $
   | b^j(x) - b(x_j) | < \varepsilon,
   $
 for any $x \in \R^d$.   Let us consider
 $$
 {\mathcal L}_j \varphi (x) =  {\rm p. v.}\int_{\R^d\backslash\{0\}
} [\varphi(x+  z)-\varphi(x)] \nu (dz)+\langle b^j(x),D\varphi(x)\rangle.  
 $$  
 Such operators verifies \A{C${}_\varepsilon $} and so by the first part of the proof,  for any initial distribution $\mu$ the  stopped martingale problem for $({\mathcal L}_j, U_j, \mu)$ is well-posed, for any $j \ge 1$. Thanks to \eqref{bo1}  the  stopped martingale problem for $({\mathcal L}, U_j, \mu)$ is also well-posed, for any $j \ge 1$.
    This finishes the proof.

 \appendix
\section{Proof of Lemma \ref{LEMME_DENS}}
We recall that we here aim at controlling the density and its derivatives of the random variable
$Z_t$ where $Z$ is a stable process of index $\alpha=1 $ satisfying the non-degeneracy condition \eqref{ND_MES} in assumption \A{ND}.


\smallskip  
 Let us recall that, for a given fixed $t>0$, we can use an It\^o-L\'evy  decomposition
 at the associated  characteristic stable time scale (i.e. the truncation is performed at the threshold $t $) 
to write $Z_t:=M_t^t+N_t^t$
where $M_t^t$ and $N_t^t $ are independent random variables. 
More precisely, 
 \begin{equation} \label{dec}
 N_s^t = \int_0^s \int_{ |x| > t }
\; x  P(du,dx), \;\;\; \; M_s^t = Z_s - N_s^t, \;\; s \ge 0,
 \end{equation} 
where $P$ is the  Poisson random measure associated with the process $Z$; for the considered fixed $t>0$,
 $M_t^t$ and $N_t^t$ correspond to
 the \textit{small jumps part } and
\textit{large jumps part} respectively w.r.t. the corresponding \textit{typical scale} of order $t$. 
A similar decomposition has been already used in the literature 
(see, for instance  the proof of  Lemma 4.3 in \cite{huan:meno:prio:19} and the references therein).
It is useful to note that the cutting threshold in \eqref{dec} precisely yields for the considered $t>0$ that:
\begin{equation} \label{ind}
N_t^t  \overset{({\rm law})}{=} t N_1^1 \;\; \text{and} \;\;
M_t^t  \overset{({\rm law})}{=} t M_1^1.
\end{equation}  
To check the assertion about $N^t$ we start with 
$$
\E [e^{i \langle p , N_t^t \rangle}] = 
\exp \Big(  t
\int_{\S^{d-1}} \int_{t}^{\infty}
 \Big(\cos (\langle p, r\theta \rangle)  - 1  \Big) \, \frac{dr}{r^{2}}\tilde \mu(d\theta) \Big), \;\; p \in \R^d
$$
(see  \eqref{GEN_DRIVING_Z} and \cite{sato:99}). Changing variable to $\frac{r}{t} =s$
we get that $\E [e^{i \langle p , N_t^t \rangle}]$ $= \E [e^{i \langle p , t N_1^1 \rangle}]$ for any $p \in \R^d$ and this shows the assertion (similarly we get
the statement for $M$).
The density of $Z_t$ then writes
\begin{equation}
\label{DECOMP_G_P}
p_{Z_t}(x)=\int_{\R^d} p_{M^t_t}(x-\xi)P_{N_t^t}(d\xi),
\end{equation}
where $p_{M^t_t}(\cdot)$ corresponds to the density of $M_t^t$ and $P_{N_t^t}$ stands for the law of $N_t^t$. 
{From Lemma A.2 in \cite{huan:meno:prio:19} (see as well
Lemma B.1  }
in \cite{huan:meno:15}),   $p_{M^t_t}  (\cdot)$  belongs to the Schwartz class ${\mathscr S}(\R^d) $ and satisfies that for all $m\ge 1 $ and all 
multi-index  \textbf {$   \beta$}, $|{\bf \beta}|\le 2 $,  there exist constants  $\bar C_m,\ C_{m}$ s.t. for all $t>0,\ x\in  \R^d  $:
\begin{equation}
\label{CTR_DER_M}
|D_x^\beta p_{M^t_t}(x)|\le \frac{\bar C_{m}}{t^{\ell }} \, p_{\bar M}(t,x),\;\; \text{where} \;\; p_{\bar M}(t,x)
:=
\frac{C_{m}}{t^{d}} \left( 1+ \frac{|x|}{t}\right)^{-m},
\end{equation}
where $C_m$ is chosen in order that {\it $p_{\bar M}(t,\cdot ) $ be a probability density.}

We carefully point out that, to establish the indicated results, since we are led to consider potentially singular spherical measures,  we only focus on integrability properties similarly to \cite{huan:meno:prio:19}. The main idea thus consists in exploiting 
{\eqref{dec},}  \eqref{DECOMP_G_P} and \eqref{CTR_DER_M}.
The derivatives on which we want to obtain quantitative bounds  will be expressed through derivatives of $p_{M^t_t}(\cdot)$, which also give the corresponding time singularities. However, as for general stable processes, the integrability restrictions come from the large jumps (here $N_t^t $) and only depend on its stability index here equal to 1.
A crucial point then consists in observing  that the convolution $\int_{\R^d}p_{\bar M}(t,x-\xi)P_{N_t^t}(d\xi) $ actually corresponds to the density of the random variable 
\begin {equation} \label{we2}
\bar Z_t:=\bar M_t+N_t^t,\;\; t>0 
\end{equation}
 (where $\bar M_t $ has density $p_{\bar M}(t,.)$ and is independent of $N_t^t $; 
 {to have such decomposition one can define each $\bar Z_t$ on a product probability space}). Then, the integrability properties of $\bar M_t+N_t^t $, and more generally of all random variables appearing below, come from those of $\bar M_t $ and $N_t^t$. 
 
 {\sl The function $\bar q(t,\cdot)$ will be  the density of the random variable $\bar Z_t$, $t>0$. }

It is readily seen that $p_{\bar M}(t,x) = {t^{- d}} \, p_{\bar M} (1, t^{- 1} x),$ $ t>0, \, $ $x \in \R^d.$  Hence 
$$
\bar M_t  \overset{({\rm law})}{=} t \bar M_1,\;\;\; N_t^t  \overset{({\rm law})}{=} t  N_1^1.
$$
By independence of $\bar M_t$ and $N_t^t$, using the Fourier transform, one can  prove that 
\begin{equation} \label{ser1}
\bar Z_t  \overset{({\rm law})}{=} t \bar Z_1.
\end{equation} 
Moreover, 
$
\E[|\bar Z_t|^\gamma]=\E[|\bar M_t+N_t|^\gamma]\le C_\gamma t^{\gamma }(\E[|\bar M_1|^\gamma]+\E[| N_1^1|^\gamma])\le C_\gamma t^{\gamma }, \; \gamma \in (0,1).
$ 
This shows that the density of $\bar Z_t$ verifies \eqref{INT_BOUND}. 

The controls \eqref{DER_DENS_BOUND} on the derivatives are derived similarly using 
\eqref{CTR_DER_M} for  all 
multi-index  \textbf {$   \beta$}, $|{\bf \beta}|\le 2 $, and the same previous argument.

Now, the bounds of \eqref{LAP_FRAZ} involving diagonal perturbation again follow from the expression of $p_{\bar M}(t,\cdot) $ in \eqref{CTR_DER_M}. Similar arguments apply to get \eqref{PERTURB_DIAG_DENS_BD}. Also, still in \eqref{LAP_FRAZ}, the bound on the fractional Laplacian is a consequence of the previous decomposition applying the operator  $\Delta^{\frac 12} $ to  $p_{M^t_t}(\cdot)$. Namely, it is easily checked that $|\Delta^{\frac 12} p_{M^t_t}(x)|\le Ct^{-1}p_{\bar M}(t,x) $ (see again Lemma 4.3 in \cite{huan:meno:prio:19} for details). Equations \eqref{DIFF_LAP_GRAD} and \eqref{DIFF_LAP_FRAZ} eventually follow from the previous bounds introducing the diagonal/off-diagonal cut-off. Namely, for \eqref{DIFF_LAP_GRAD}, if $|z-z'|>t $, then  $|Dp_{Z_t}(z)-Dp_{Z_t}(z')|\le |Dp_{Z_t}(z)|+|Dp_{Z_t}(z')|\le \frac{C}t\Big(\frac{|z-z'|}{t} \Big)^\beta \big(\bar q(t,z)+ \bar q(t,z')\big)  $, whereas if $ |z-z'|\le t$ then, from \eqref{DER_DENS_BOUND} $|Dp_{Z_t}(z)-Dp_{Z_t}(z')| \le \int_0^1 \bar q(t,z+\lambda (z'-z)) \frac{|z-z'|}{t^2}\le \frac Ct\bar q(t,z) \big(\frac{|z-z'|}{t}\big)^\beta$, using \eqref{PERTURB_DIAG_DENS_BD} for the last inequality. The bound \eqref{DIFF_LAP_FRAZ} can be derived in a similar way. \qed

  \begin{lemma} \label{d66}
   Let  $b \in C^2_b(\R^d, \R^d)$ and assume  that there exists $b_0 \in \R^d$
  $\varepsilon \in (0,1)  $, such that, 
for all $x\in \R^d$,
$ 
|b(x)-b_0 | \le \varepsilon
$ (cf. \A{C${}_\varepsilon$}). 
    Let us consider the pathwise unique solution $(X_t^x)$ to 
  \begin{gather*}
X_t = x + \int_0^t b(X_s) ds + Z_t, \;\; t \ge 0
\end{gather*}
(defined on a stochastic basis $(\Omega, {\mathcal  F}, ( {\mathcal  F}_t), \P)$) and the corresponding  resolvent 
$$
u_{\lambda}(x) = {{\E}}\Big[\int_0^{+\infty}  \exp(-\lambda s )  f(X_s^x) ds\Big], \;\; \lambda >0,\; x \in \R^d.
$$
Then, for 
$\lambda \ge 1$ and $p>d$, there exists $C = C( \varepsilon, d, p, \| b \|_{C^2_b})$ such that 
\begin{equation}\label{dd}
|u_{\lambda}(x)| \le C \| f\|_{L^p({\R}^d)}, \;\; f \in C_0^{\infty}(\R^d).
\end{equation} 
  \end{lemma}
\begin{proof}
 Thanks to the regularity of $b$ we know that $u =u_{\lambda}$ is the unique bounded classical solution to 
 \begin{gather*}
 \lambda u(x) - L(x) - b(x) \cdot Du(x) = f(x),\;\; x \in \R^d,
\end{gather*}
which we can write as $\lambda u(x) - L(x) - b_0 \cdot Du(x) = f(x) + (b(x)-b_0) \cdot Du(x)$. It follows the representation formula
\begin{equation*}\label{d}
u(x)=  \int_0^{+\infty}  \exp(-\lambda t ) dt\int_{\R^d} [f(y) + (b(y)-b_0) 
\cdot Du(y)] \, p_{Z_t} (y- x - tb_0) dy. 
\end{equation*}
Now we use 
$ |b(x)-b_0 | \le \varepsilon$ for $\varepsilon$ small enough and   \eqref{BOUND_LP_GRAD}.  By a fixed point theorem in $W^{1,p}(\R^d)$, $1< p< \infty$, we find that $u \in W^{1,p}(\R^d)$ and, moreover, 
\begin{equation}\label{dd}
\|u \|_{W^{1,p}(\R^d)} \le C \| f\|_{L^p(\R^d)}.
\end{equation} 
Choosing $p>d$ and applying the Sobolev embedding theorem we obtain the assertion.
\end{proof}

\end{document}